\documentclass[12pt]{article}
\usepackage{amsmath,latexsym,amssymb}
\usepackage{eucal}
\newtheorem{thm}{Theorem}[section]
\newtheorem{lem}[thm]{Lemma}

\newtheorem{df}[thm]{Definition}
\newtheorem{cor}[thm]{Corollary}
\newtheorem{exam}[thm]{Example}
\newcommand{\id}{\mathrm{id}}

\newcommand{\rAC}{\mathrm{AC}}
\newcommand{\rB}{\mathrm{B}}

\newcommand{\cM}{\mathcal{M}}
\newcommand{\cN}{\mathcal{N}}

\newcommand{\cQ}{\mathcal{Q}}
\newcommand{\cR}{\mathcal{R}}

\newcommand{\tM}{\tilde{\mathcal{M}}}
\newcommand{\tN}{\tilde{\mathcal{N}}}
\newcommand{\tQ}{\tilde{\cQ}}

\newcommand{\fH}{\mathfrak{H}}

\DeclareMathOperator{\Ad}{Ad}

\DeclareMathOperator{\Int}{Int}

\DeclareMathOperator{\End}{End}

\begin{document}

\title{On the relative bicentralizer flows  and the relative flow of
weights of inclusions of factors of type III$_1$}
\author{Toshihiko MASUDA\footnote{Supported by 
JSPS KAKENHI Grant Number 16K05180.} \\
 Graduate School of Mathematics, Kyushu University \\
 744, Motooka, Nishi-ku, 
Fukuoka, 819-0395, JAPAN \\
e-mail address: masuda@math.kyushu-u.ac.jp}
\date{}
\maketitle

\begin{abstract}
 We show the relative bicentralizer flow and the relative flow of
 weights are isomorphic for an inclusion of injective factors of type
 III$_1$ with finite index, or an irreducible discrete inclusion whose
 small algebra is an injective factor of type III$_1$.
\end{abstract}

\section{Introduction}\label{sec:intro}

It is well known that the uniqueness of the injective factor of type
III$_1$ is solved by \cite{Co-III1}, \cite{Ha-III1}, that is, 
A. Connes tried to show the uniqueness of the injective factor
of type III$_1$, and found that an injective factor of type III$_1$ with 
trivial bicentralizer algebra is isomorphic to the Araki-Woods factor of
type III$_1$ \cite{AW}. Then U. Haagerup \cite{Ha-III1} solved the
bicentralizer problem, and the uniqueness of the injective factor of
type III$_1$ is established. 
It is natural to ask if 
a bicentralizer algebra of  any type III$_1$ factor is trivial.
Although this is still an open problem, no counterexample has been known
so far.
(See the introduction of \cite{AHHM-bicent} and reference 
therein about recent progress on the trivial bicentralizer problem.)

In \cite{AHHM-bicent}, H. Ando, U. Haagerup, C. Houdayer and 
A. Marrakchi showed that one can associate a flow on a (relative) bicentralizer
algebra canonically, and investigated its properties. On the other hand,
there exists a canonical flow, i.e., the trace scaling action on the
continuous core. Thus one can define the notion of a relative flow of
weights for an inclusion of factors. (This notion was appeared in
\cite{Win-flow} related to the  classification of subfactors.)
They ask whether these two flows are conjugate, and
present an example with positive answer. 

In this paper, we give an affirmative answer when a small algebra of a
given inclusion is an injective factor of type III$_1$. Our proof is
based on Popa's solution of relative version of bicentralizer problem
\cite{Po-CBMS}. 

\textbf{Acknowledgements.} The author is grateful to Professor H. Ando,
and C. Houdayer for comments on this paper. He also thanks
Prof. A. Marrakchi for pointing out  Corollary  \ref{cor:bicent} to him.

\section{Preliminaries and notations}

Throughout this paper, we assume that all von Neumann algebras have
separable predual.
For a von Neumann algebra $\cM$, we denote its continuous core by $\tM$,
the implementing unitary for $\sigma^\varphi_s$ by 
$\lambda^\varphi(s)\in \tM$, and the trace scaling automorphism by $\theta_t$.

Let $\End(\cM)$ be the set of all unital endomorphisms of $\cM$, and 
$\End_0(\cM)=\{\rho\in \End(\cM)\mid [\cM:\rho(\cM)]<\infty\}$. 
We denote $[\cM:\rho(\cM)]^{\frac{1}{2}}=d(\rho)$ for $\rho\in \End_0(\cM)$.
Define the left inverse of $\rho$ by $\phi_\rho:=\rho^{-1}E_\rho$,
where $E_\rho$ is the minimal conditional expectation $\cM\rightarrow \rho(\cM)$.
For $\rho\in \End_0(\cM)$, let $\tilde{\rho}\in \End(\tM)$  be the canonical
extension \cite{Iz-can2}. 
For $\rho,\sigma\in \End(\cM)$, the intertwiner space $(\rho,\sigma)$ is
defined by $(\rho,\sigma):=\{T\in \cM\mid \sigma(x)T=T\rho(x), x\in \cM\}$.

Let $\cN\subset \cM $ be an inclusion of von Neumann algebras with a conditional
expectation $E:\cM\rightarrow \cN$. When $\cN\subset \cM$ is a subfactor
with finite index,  we always choose $E$ as a minimal conditional expectation. 
We recall the notion of a relative bicentralizer \cite{Ha-III1},
\cite{M-III1}.
(We use the notations in \cite{AHHM-bicent}.)
\begin{df}\label{df:bicent}
 Let $\varphi$ be a faithful normal state on $\cM$ with $\varphi\circ
 E=\varphi$. \\
$(1)$ The asymptotic centralizer of $\varphi$ is defined by
\[
 \rAC(\cN,\varphi)=\{(x_n)\in\ell^\infty(\mathbb{N},\cN)\mid
 \lim_{n\rightarrow \infty}\|x_n\varphi-\varphi x_n\|=0\}
\]
$(2)$ The relative bicentralizer algebra of $\varphi$ is defined by
\[
\rB(\cN\subset \cM,\varphi)=\left\{a\in \cM\mid 
\lim\limits_{n\rightarrow \infty}[x_n,a]=0 \mbox{ strongly for
 all } (x_n)\in \rAC(\cN,\varphi). \right\}
\]
When $\cN=\cM$, we simply denote it by $\rB(\cM,\varphi)$.
\end{df}

The following theorem due to S. Popa is essential in our argument, whose proof is an
appropriate modification of that of Haagerup \cite{Ha-III1}.
(Also see \cite[Corollary 4.7]{M-III1}.)
\begin{thm}[{\cite[Theorem 4.2]{Po-CBMS}}]\label{thm:Popabicent}
 Let $\cN\subset \cM$ be an inclusion of injective factors of type
 III$_1$ with finite index. Assume that $\cN'\cap \cM=\tN'\cap
 \tM$. Then $\rB(\cN\subset \cM,\varphi)=\cN'\cap \cM$.
\end{thm}


Let $\varphi,\psi$ be faithful normal states on $\cM$ with
$\varphi=\varphi\circ E$, $\psi=\psi\circ E$.
In \cite[Theorem A]{AHHM-bicent}, the authors showed the existence of 
a flow $\beta:\mathbb{R}_+^*\curvearrowright\rB(\cN\subset
\cM,\varphi)$, and the canonical isomorphism $
\beta^{\psi,\varphi}:B(\cN\subset \cM,\varphi)\rightarrow B(\cN\subset
\cM,\psi)$, 
which intertwines $\beta^\varphi$ and $\beta^\psi$. 
We call $\beta_\lambda^\varphi$ a relative bicentralizer flow.

The flow  $\beta_\lambda^\varphi$ 
is characterized as follows;
$\lim\limits_{n\rightarrow \infty}a_nx-\beta_\lambda^\varphi(x)a_n=0 $ in the strong topology for
all $(a_n)\in \ell^\infty (\mathbb{N},\cN)$ with
$\lim\limits_{n\rightarrow \infty}\|a_n\varphi-\lambda \varphi a_n\|=0$.

In a similar way
the isomorphism $\beta^{\psi,\varphi}$ is characterized as follows;
$\lim\limits_{n\rightarrow \infty}a_nx-\beta^{\psi,\varphi}(x)a_n=0 $ in the strong topology for
all $(a_n)\in \ell^\infty (\mathbb{N},\cN)$ with $\lim\limits_{n\rightarrow \infty}\|a_n\varphi-\psi a_n\|=0$.

\section{Main results}\label{sec:bicent}
We begin with the following lemma, which is key for our main results.
\begin{lem}\label{lem:bicent}
 Let $\cN$ be an injective factor of type III$_1$, and 
 fix a faithful normal state $\varphi$ on $\cN$.
Assume that $\rho\in \End_0(\cN)$ is a non-modular endomorphism in the sense
 of \cite{Iz-can2}, and 
$a\in \cN$ satisfies $x_na-a\rho(x_n)\rightarrow 0$
 in the strong topology for all $\{x_n\}\in
 \rAC(\cN,\varphi)$. Then $a=0$. 
\end{lem}
\textbf{Proof.} Set $\cQ:=\cN\otimes M_2(\mathbb{C})$, and define an embedding
$\pi:\cN\rightarrow \cQ$ by $\pi(x)=x\otimes e_{11}+\rho(x)\otimes
e_{22}$. 
Define a faithful normal state $\Psi$ on $\cQ$   and a conditional
expectation $E_\Psi:\cQ\rightarrow \pi(\cN)$ by 
\[
 \Psi((a_{ij}))=\frac{1}{1+d(\rho)}\left(a_{11}+d(\rho)\phi_\rho(a_{22})
 \right),\,\,
E_\Psi((a_{ij}))=\frac{1}{1+d(\rho)}\pi(a_{11}+d(\rho)\phi_\rho(a_{22})).
\]
Then $E_\Psi$ is the minimal expectation with $\mathrm{Ind}(E_\Psi)=1+d(\rho)$, and $\Psi=\Psi\circ E_\Psi$.
It is easy to see $\varphi(x)=\Psi(\pi(x))$ for $x\in \cN$. Thus we have 
$\pi\left(\rAC(\cN,\varphi)\right)=\rAC(\pi(\cN),\Psi)$. 

The core inclusion for $\pi(\cN)\subset  \cQ$ is given by
$\tilde{\pi}(\tN)\subset \tQ=\tN\otimes M_2(\mathbb{C})$, where 
$\tilde{\pi}(x)=x\otimes e_{11}+\tilde{\rho}(x)\otimes  e_{22}$. 
The irreducibility and non-modularity of $\rho$ implies
$(\id,\tilde{\rho})=\{0\}$ \cite[Proposition 3.4(1)]{Iz-can2}, and yields
\[
 \pi(\cN)'\cap \cQ=\tilde{\pi}(\tN)'\cap \tQ=\mathbb{C}e_{11}\oplus \mathbb{C}e_{22}.
\]
By Theorem \ref{thm:Popabicent},
\[
 \rB(\pi(\cN)\subset \cQ,\Psi)=
 \pi(\cN)'\cap \cQ=\mathbb{C}e_{11}\oplus \mathbb{C}e_{22}.
\]

Assume $a\in \cN$ satisfies the condition in the statement of Lemma. Then 
$a\otimes e_{12}\in \rB(\pi(\cN)\subset \cQ,\Psi)$, and hence
       $a=0$. \hfill$\Box$

\bigskip

\noindent
\textbf{Remark.} For a type III$_1$ factor $\cN$, $\rho\in \End_0(\cN)$
is an irreducible modular endomorphism if and only if $\rho=\Ad u\circ
\sigma^\varphi_t$ for some $u\in U(\cN)$. \\

\medskip

Let $\cN$ be an injective factor of type III$_1$, and 
$\cN\subset \cM$ a (possibly reducible) inclusion of factors of type
III$_1$ with finite index, or a discrete irreducible inclusion. 
 Fix a faithful normal state on $\cM$ such
that $\varphi=\varphi\circ E$. 
We recall results in \cite{ILP-Galois} on discrete inclusions.
Let $\gamma|_{\cN}=\bigoplus_{\xi \in \Xi}
n_\xi[\rho_\xi ]$ be an irreducible decomposition of the canonical
endomorphism $\gamma$ for $\cN\subset \cM$, and $\Xi_0:=\{t\in \mathbb{R}\mid \sigma^\varphi_t \prec
\gamma|_{\cN} \}$. We regard $\Xi_0$ as a subset of $\Xi$. 
We fix a representative of $[\rho_\xi]$ so that
$\rho_s=\sigma^\varphi_s$ for $s\in \Xi_0$. 
Let 
\[
 \fH_\xi:=\{V\in \cM \mid Vx=\rho_\xi(x)V, \forall x\in \cN\}.
\]
Then $\fH_\xi$ is a finite dimensional Hilbert space with $\dim \fH_\xi=n_\xi$
by the inner product $\langle V,W\rangle=E(W^*V)$. 
Fix a CONS $\{V^\xi_k\}_{k=1}^{n_\xi}\subset \fH_\xi$.
Then $a\in \cM $ can be expand as $a=\sum_{\xi,k}a_{k,\xi}V^\xi_k$,
$a_{k,\xi}\in \cN$,  uniquely. 
When $\cN\subset \cM$ is of
infinite index, this expansion does not converge in the usual
operator algebra topology. However coefficients $\{a_{k,\xi}\}_{k,\xi}$ determines $a$ uniquely.

\begin{thm}\label{thm:bicent}
 Let $\cN\subset \cM$ be as above. \\
$(1)$ The relative bicentralizer algebra $B(\cN\subset \cM,\varphi)$ is given by
\[
 \rB(\cN\subset \cM,\varphi)=\left\{\sum_{s \in \Xi_0}\sum_{k=1}^{n_s} a_{k, s }V^s_{k}\mid a_{k,s}\in \mathbb{C}\right\}.
\]
The restriction of $E$ on $\rB(\cN\subset \cM,\varphi)$ gives a faithful
 normal tracial state. \\
$(2)$ The relative bicentralizer flow $\beta^\varphi_\lambda$ is given by $
\beta^\varphi_\lambda(V)=\lambda^{is}V$, $V\in \fH_s$. \\
$(3)$ Let $\psi$ be another faithful normal state on $\cM$ with
 $\psi=\psi\circ E$. 
The isomorphism $\beta^{\psi,\varphi}:B(\cN\subset
 \cM,\varphi)\rightarrow B(\cN\subset \cM,\psi)$
is given by
$\beta_{\psi,\varphi}(V)=[D\psi:D\varphi]_sV$ for $V\in \fH_s$.
\end{thm}
\textbf{Proof.} (1) Take $a\in \rB(\cN\subset \cM,\varphi)$. Thus
$x_na-ax_n\rightarrow 0$ in the strong topology for all
$\{x_n\}\in \rAC(\cN,\varphi)$.  By expanding $a=\sum_{k,\xi}a_{k,\xi}V^\xi_k$,
we get $x_na_{k,\xi}-a_{k,\xi}\rho_{\xi}(x_n)\rightarrow 0$ in the
strong topology. By Lemma
\ref{lem:bicent}, $a_{k,\xi}=0$ for $\xi\not\in \Xi_0$. If $s\in
\Xi_0$, then we can easily see $a_{k,s}\in \mathbb{C}$ by Haagerup's
theorem $\rB(\cN,\varphi)=\mathbb{C}$ \cite{Ha-III1} and the fact
$\lim\limits_{n\rightarrow \infty}\sigma_s^\varphi(x_n)-x_n=0$ in the
strong topology. 

Take $a=\sum\limits_{s \in \Xi_0}\sum\limits_{k=1}^{n_s} a_{k, s }V^s_{k}\in
\rB(\cN\subset \cM,\varphi)$. 
We choose a CONS $\{V_k^0\}_{k=1}^{n_0}\subset \fH_0$ so that $V_1^0=1$. (If $\cN\subset
\cM$ is irreducible, $\dim\fH_0=1$.)
Then $E(a)=a_{1,0}$, and we get
$E(a^*a)=E(aa^*)=\sum\limits_{s\in \Xi_0}\sum\limits_{k=1}^{n_s}|a_{k,s}|^2$. Hence
the restriction of $E$ on $\rB(\cN\subset \cM,\varphi)$ is tracial.

(2) 
The proof is similar to that of \cite[Proposition 6.11]{AHHM-bicent}.
Take $(a_n)\in \ell^\infty (\mathbb{N},\cN)$ with
$\lim_{n}\|a_n\varphi-\lambda\varphi a_n\|=0$. Then
$\sigma^\varphi_{s}(a_n)-\lambda^{-is}a_n\rightarrow 0$, and hence we have
\begin{align*}
a_n V-\lambda^{is}Va_n&= \lambda^{is}\left(
\lambda^{-is}a_n -\sigma^\varphi_{s}(a_n)\right)V\rightarrow 0,\,\, V\in \fH_s.
\end{align*}
This shows $\beta^\varphi_\lambda(V_k^s)=\lambda^{is}V_k^s$. 

(3)
Take $(a_n)\in \ell^\infty(\mathbb{N},\cN)$ with $\|a_n\varphi-\psi
a_n\|\rightarrow 0$. Then
$[D\psi:D\varphi]_s\sigma^\varphi_s(a_n)-a_n\rightarrow 0$ in the
strong topology.
(To see this, consider $\cN\otimes M_2(\mathbb{C})$ and a state
$\Psi((a_{i,j}))=\dfrac{1}{2}(\varphi(a_{11})+\psi(a_{22}))$.)
For $V\in \fH_s$, we have
\begin{align*}
a_nV-[D\psi:D\varphi]_sV
 a_n=(a_n-[D\psi:D\varphi]_s\sigma^\varphi_s(a_n))V\rightarrow 0.
\end{align*}
This means $\beta_{\psi,\varphi}(V)=[D\psi:D\varphi]_sV$ for $V\in \fH_s$.
\hfill$\Box$

\medskip

\noindent
\textbf{Remark.} Even if $\cN\subset \cM$ is of infinite index and
reducible, the same proof works if ``the Fourier expansion'' holds, e.g.,
$\cN\subset \cM$ is a crossed product by a discrete group action.

\bigskip

Combining Theorem \ref{thm:bicent} and \cite[Theorem 3.5]{HI-bicent}, we get the
following.
(This is pointed out by A. Marrakchi to the author.)
\begin{cor}\label{cor:bicent}
 Let $\cN\subset \cM$ be 
 an irreducible discrete inclusion of factors of type III$_1$ such that
 $\cN$ is injective. Then the bicentralizer algebra
 $\rB(\cM,\varphi)$ of $\cM$ is trivial. 
\end{cor}
\textbf{Proof.} By definition, it is trivial $\rB(\cM,\varphi)\subset
\rB(\cN\subset \cM,\varphi)$. Since $\rB(\cN\subset \cM,\varphi)$ is finite,
\cite[Theorem 6.5]{HI-bicent} implies $\rB(\cM,\varphi)=\mathbb{C}$. \hfill$\Box$
 
\bigskip

In \cite{AHHM-bicent}, the authors 
ask whether the relative bicentralizer flow 
is conjugate to the relative flow of weights $\{\cN'\cap \tM,
\theta_t\}$, and they exhibit examples in \cite[Proposition
6.11]{AHHM-bicent}. Our case also provides a positive answer to their
question.

\begin{thm}\label{thm:bicentflow} Let $\cN\subset \cM$ be as in Theorem \ref{thm:bicent}.
Then two flows 
 $\{\rB(\cN\subset \cM,\varphi),\beta^\varphi_\lambda\}$ and  $\{\cN'\cap \tM, \theta_t\}$ are
 conjugate via a map $\alpha$ given by
 $\alpha(V)=\lambda^{\varphi}(s)^*V$, $s\in \Xi_0$, $V\in\fH_s$.
\end{thm}
\textbf{Proof.}
At first we treat the  case $[\cM:\cN]<\infty$, and show the following
claim. 

\medskip

\noindent
\textbf{Claim.} 
When $\cN\subset \cM$ is of finite index, 
we have $\sigma^\varphi_t(V)=d(\rho_\xi)^{-it}[D\varphi\circ \phi_\xi:D\varphi]^*_tV$ for any $V\in \fH_\xi$. 
where $\phi_\xi$ is the left inverse of $\rho_\xi$.
In particular, $\sigma^\varphi_t(V)=V$ for $V\in \fH_s$, $s\in \Xi_0$.
\medskip

This statement is proved in \cite[p.45]{ILP-Galois}. We present
a proof for readers convenience.

We recall the
following general fact. For injective
unital  homomorphisms $\rho,\sigma:\cN\rightarrow \cM$ with finite
index, take $T\in \cM$ satisfying $\rho(x)T=T\sigma(x)$, $x\in \cN$.
Then we have 
\[
 d(\rho)^{it}[D(\psi\circ \phi_\rho):D\chi]_t
\sigma^{\chi}_t(T)=T
 d(\sigma)^{it}[D(\psi\circ \phi_\sigma):D\chi]_t
\] for faithful normal states $\psi\in \cN_*$, $\chi\in \cM_*$, cf. \cite[Lemma A.2]{MaTo-JFA}.
(We can define $d(\rho)$ and $\phi_\rho$ for $\rho:\cN\rightarrow \cM$ similarly.)

Let $\iota:\cN\rightarrow \cM$ an inclusion map, and apply the above
result for $\rho=\iota\rho_\xi$, $\sigma=\iota$, $T=V\in \fH_\xi$,
$\psi=\varphi|_{\cN}$, $\chi=\varphi$.
Since left inverses are given by 
$\phi_{\iota\rho_\xi}=\phi_\xi \circ E$,
$\phi_{\iota}= E$, we get the desired result.

Claim implies that $Va=\tilde{\rho_\xi}(a)V$ for $a\in \tN$,
$V\in \fH_\xi$. Hence $\{V_k^\xi\}_{\xi,k}$ is a basis for $\tN\subset \tM$, and
every $a\in \tM$ can be expanded $a=\sum_{k,\xi}a_{k,\xi}V_k^\xi$
uniquely. Then we get 
\[
\cN'\cap \tM=\left\{\sum_{s\in
\Xi_0}\sum_{k=1}^{n_s}a_{k, s}\lambda^{\varphi}(s)^*V_k^s\mid
a_{k,s}\in\mathbb{C}\right\} 
\]
by Connes-Takesaki relative commutant theorem \cite{CT}
and \cite[Proposition 3.4(1)]{Iz-can2}. (This also provides a
proof of the relative commutant theorem \cite[Theorem 4.3]{Po-CBMS} for
subfactors with finite index.)

Now we  can show $\alpha$ is a $*$-isomorphism. 
Take $V\in \fH_s$ and $W\in \fH_t$
for $s,t\in \Xi_0$.  On one hand, we have $VW\in \fH_{s+t}$ and $V^*\in \fH_{-s}$.
(Hence $VW=0$ if $s+t\not\in \Xi_0$.)
On the other hand, we have
\[
 \lambda^\varphi(s)^*V\lambda^\varphi(t)^*W=
 \lambda^\varphi(s)^*\lambda^\varphi(t)^*\sigma^\varphi_t(V)W=
 \lambda^\varphi(s+t)^*VW
\]
and 
\[
 (\lambda^\varphi(s)^{*}V)^*=V^*\lambda^\varphi(s)=\lambda^\varphi(s)\sigma^\varphi_{-s}(V^*)=
\lambda^{\varphi}(-s)^*V^*
\]
by Claim.
This shows $\alpha$ is a $*$-isomorphism.

Since $\theta_t(\lambda^\varphi(s)^*V_k^s)=e^{ist}\lambda^\varphi(s)^*V_k^s$, 
$\alpha\circ \beta_\lambda^\varphi\circ \alpha^{-1}=\theta_{\log
\lambda}$ holds, and hence $\alpha$ intertwines two flows $\beta^\varphi$ and $\theta$.

Let us assume that $\cN\subset \cM$ is an irreducible discrete
inclusion with $[\cM:\cN]=\infty$.
One should note that the statement of Claim may fail in this case, i.e.,
 the operator $a_\xi$ defined in \cite[p.39]{ILP-Galois} is not trivial
 in general.
(See \cite[p.45, Appendix]{ILP-Galois} for such an example.
Crossed product inclusions
by a discrete quantum group of non-Kac type also provide such examples.) 
However, one has $a_s=1$ and 
the statement of Claim holds for $s \in \Xi_0$, since $d(\rho_s)=d(\sigma^\varphi_s)=1$.
Moreover we have $\dim \fH_s=1$, $s\in \Xi_0$.

In this case, $\cN'\cap \tM$ is described in \cite[Theorem 6.9]{AHHM-bicent}, which
can be regard as a generalization of computation in \cite{Kw-Tak}, \cite{Se-flow}.
Namely, fix a unitary $V^s\in \fH$. Then we have
\[
\cN'\cap \tM=\left\{\sum_{s\in \Xi_0}a_{s}\lambda^{\varphi}(s)^*V^s\mid
a_{s}\in\mathbb{C}\right\}. 
\]
Thus the proof of finite index case also works in this case.

One should note that the restriction of a conditional expectation on
$\rB(\cN\subset \cM,\varphi)$ 
(resp. on $\cN'\cap \tM$) gives a faithful normal tracial state, and $\alpha$
preserves these states.
\hfill$\Box$

\bigskip

We present two examples of inclusions $\cN\subset \cM$ such that $\cN$ is an
injective factor of type III$_1$, but $\cM$ is not injective, and 
$\cN'\cap \cM\ne B(\cN\subset \cM,\varphi)$.
\begin{exam}
\upshape
This example is a modification of
 \cite[Proposition 6.11]{AHHM-bicent}.
 Let $N$ be an abelian discrete group such that there exist an
 injective homomorphism $\nu:N\rightarrow \mathbf{R}$, and a 2-cocycle
 $\mu\in Z^2(N,\mathbf{T})$ such that $\mu(g,h)\overline{\mu(h,g)}$ is a
 non-degenerate bicharactor. (For example, 
 $N=\mathbf{Z}\oplus  \mathbf{Z}$, $\nu(a,b)=a+\sqrt{2}b$, 
$\mu((a,b),(c,d))=e^{2\pi i \theta bc}$, $\theta\in \mathbb{R}\backslash \mathbb{Q}$.)

Let $K$ be a non-amenable discrete group,
 and $\gamma$ a free action of $K$ on an injective factor of type
 II$_1$ $\cR_0$. Then $\mu$ can be regard as an element of $Z^2(N\times
 K,\mathbb{T})$ naturally.

Let us define an action of $\alpha$ of $N\times K$ on $\cN\otimes \cR_0$ by
 $\alpha_{(n,g)}=\sigma^\varphi_{\nu(n)}\otimes \gamma_g$, and 
$\cM:=(\cN\otimes \cR_0)\rtimes_{\alpha,\mu}(N\times K)$. 
Since $\alpha$ is a free action of a non-amenable discrete group, $\cM$
 is not injective and $\cN'\cap \cM=\mathbb{C}$.
Let $V_g\in \cM$ be the implementing unitary.
By Theorem \ref{thm:bicent}, $B(\cN\subset
 \cM,\varphi)=\{\sum_{n\in N}a_nV_n\mid a_n\in \mathbb{C}\}$,  
and this algebra is a factor of type II$_1$ by the choice of $\mu$.

Since the canonical extension $\tilde{\alpha}_g$ is inner if and only if
 $g\in N$,
we can see that 
$\cN'\cap \tM=\{\sum_{n\in N}a_n\lambda^\varphi(\nu(n))^*V_n\mid n\in
 N\}$ as in \cite{Kw-Tak}, \cite{Se-flow}, \cite[Proposition  6.11]{AHHM-bicent}. 
\end{exam}

\begin{exam}
\upshape 
Let $K$ be a non-amenable discrete group, and $\Gamma$ a wreath product
of $\mathbb{Z}$ and $K$. Namely, let $N=\bigoplus_{g\in K}\mathbb{Z}$,
and define an action of $K$ on $N$ by shift. Then $\Gamma$ is defined by
$\Gamma=N\rtimes K$. 

Let $\nu:N\ni(n(g))_{g\in K} \mapsto  \sum_{g\in K}n(g)\in \mathbb{R}$.
(Here we write the group operation of $N$ additively, but will write it
 multiplicatively in what follows.)
Then $\nu$ is a homomorphism with $\nu(g^{-1}ng)=\nu(n)$ for $g\in \Gamma$.

We can construct an action $\alpha$ of $\Gamma$ on $\cN$ with following properties. \\
(1) Let $\tilde{\alpha}_g$ be the canonical extension of $\alpha_g$.
 Then $\{g\in \Gamma\mid \tilde{\alpha}_g\in \Int(\tN)\}=N$. \\
(2) There exists a unitary $u_n\in \tN$, $n\in N$, such that 
$\tilde{\alpha}_n=\Ad u_n$, $u_nu_m=u_{nm}$, and $\tilde{\alpha}_g(u_{g^{-1}ng})=u_n$, $g\in \Gamma$. \\ 
(3) $\theta_t(u_n)=e^{-it\nu(n)}u_n$. \\
(See \cite[Proposition 22]{KwST}, \cite[Section 6.1]{M-unif-Crelle} for
existence of such actions. 
Note that what we actually need is the
existence of a free action of $\Gamma$ on the injective factor of type
II$_1$ and the amenability of $N$ in the construction, and 
$\Gamma$ does not need to be amenable.)

Let $\cM=\cN\rtimes_\alpha \Gamma$, and $V_g\in \cM$ be the
 implementing unitary.
 Property (3) means that
$v_n:=u_n\lambda^\varphi(\nu(n))^*\in\tN^\theta=\cN$, and 
$\alpha_n=\Ad v_n\circ \sigma^\varphi_{\nu(n)}$. 
By Theorem \ref{thm:bicent}, we have 
\[
 \rB(\cN\subset \cM,\varphi)=\left\{\sum_{n\in N}a_nv_n^*V_n\mid a_n\in \mathbb{C}\right\},
\]
which is isomorphic to the group von Neumann algebra of $N$.
By \cite{Kw-Tak}, \cite{Se-flow}, we can see that the center of $\tM$ is
trivial, i.e., $\cM$ is of type III$_1$. 
Indeed, we can see that 
$\cN'\cap \tM=\{\sum_{n\in N}a_nu_n^*V_n\mid a_n\in\mathbb{C}\}$
by property (1).
Then it is easy to see $\cM'\cap \tM=\mathbb{C}$ by using property (2).
Since $\Gamma$ is non amenable, $\cM$ is not injective. 

The inner part of $\alpha$ is $\mathrm{Ker}(\nu)$, and it is  not a
 trivial subgroup.
Thus we have
\[
 \cN'\cap \cM=\left\{\sum_{n\in
 \mathrm{Ker}(\nu)}a_nv_n^*V_n\mid a_n\in \mathbb{C}\right\},
\]
and hence  $\cN\subset  \cM$  is not irreducible. 
 But we can see that Theorem \ref{thm:bicentflow} holds in this case 
by the description of
 $\rB(\cN\subset \cM,\varphi)$ and $\cN'\cap \tM$.

\end{exam}


\begin{thebibliography}{10}

\bibitem{AHHM-bicent}
Ando, H., Haagerup, U., Houdayer, C., and Marrakchi, A., \textsl{Structure of
  bicentralizer algebras and inclusions of type {III} factors}, preprint,
  http://arxiv.org/abs/1804.05706,  (2018).

\bibitem{AW}
Araki, H. and Woods, J., \textsl{A classification of factors}, Publ. Res. Inst.
  Math. Sci. {\bf 3} (1968), 51--130.

\bibitem{Co-III1}
Connes, A., \textsl{Factors of type {III}$_1$, property ${L}'_\lambda$ and
  closure of inner automorphisms}, J. Operator Theory {\bf 14} (1985),
  189--211.

\bibitem{CT}
Connes, A. and Takesaki, M., \textsl{The flow of weights on factors of type
  {III}}, Tohoku J. Math. {\bf 29} (1977), 473--555.

\bibitem{Ha-III1}
Haagerup, U., \textsl{Connes' bicentralizer problem and uniqueness of the
  injective factor of type {III}$_1$}, Acta Math. {\bf 158} (1987), 95--148.

\bibitem{HI-bicent}
Houdayer, C. and Isono, Y., \textsl{Unique prime factorization and
  bicentralizer problem for a class of type {III} factors}, Adv. Math {\bf 305}
  (2017), 402--455.

\bibitem{Iz-can2}
Izumi, M., \textsl{Canonical extension of endomorphisms of type {III} factors},
  Amer. J. Math. {\bf 125} (2003), 1--56.

\bibitem{ILP-Galois}
Izumi, M., Longo, R., and Popa, S., \textsl{A {G}alois correspondence for
  compact groups of automorphisms of von {N}eumann algebras with a
  generalization to {K}ac algebras}, J. Funct. Anal. {\bf 155} (1998), 25--63.

\bibitem{KwST}
Kawahigashi, Y., Sutherland, C.~E., and Takesaki, M., \textsl{The structure of
  the automorphism group of an injective factor and the cocycle conjugacy of
  discrete abelian group actions}, Acta Math. {\bf 169} (1992), 105--130.

\bibitem{Kw-Tak}
Kawahigashi, Y. and Takesaki, M., \textsl{Compact abelian group actions on
  injective factors}, J. Funct. Anal. {\bf 105} (1992), 112--128.

\bibitem{M-III1}
Masuda, T., \textsl{An analogue of {C}onnes-{H}aagerup approach to
  classification of subfactors of type {III}$_1$}, J. Math. Soc. Japan {\bf 57}
  (2005), 959--1003.

\bibitem{M-unif-Crelle}
Masuda, T., \textsl{Unified approach to classification of actions of discrete
  amenable groups on injective factors}, J. Reine Angew. Math. {\bf 683}
  (2013), 1--47.

\bibitem{MaTo-JFA}
Masuda, T. and Tomatsu, R., \textsl{Classification of minimal actions of a
  compact {K}ac algebra with amenable dual on injective factors of type {III}},
  J. Funct. Anal. {\bf 274} (2007), 1965--2025.

\bibitem{Po-CBMS}
Popa, S., \textsl{Classification of subfactors and their endomorphisms},
  Regional {C}onference {S}eries in {M}athematics (1995), no.~86.

\bibitem{Se-flow}
Sekine, Y., \textsl{Flow of weights of the crossed products of type {III}
  factors by discrete groups}, Publ. Res. Inst. Math. Sci. {\bf 26} (1990),
  655--666.

\bibitem{Win-flow}
Winsl{\o}w, C., \textsl{The flow of weights in subfactor theory}, Publ. Res.
  Inst. Math. Sci. {\bf 31} (1995), 519--532.

\end{thebibliography}

\ifx\undefined\bysame
\newcommand{\bysame}{\leavevmode\hbox to3em{\hrulefill}\,}
\fi

\end{document}